\documentclass{amsart}
\usepackage{amsmath,amssymb, amsthm, amscd}
\input xy
\xyoption{all}

\theoremstyle{definition}

\theoremstyle{remark}

\numberwithin{equation}{section}

\begin{document}

\title{Equivariant $K$-theory and the Chern character for discrete groups}
\author{Efton Park}
\address{Department of Mathematics, Box 298900, 
Texas Christian University, Fort Worth, TX 76129}
\email{e.park@tcu.edu}
\date{June 14, 2010}

\begin{abstract}
Let $X$ be a compact Hausdorff space, let $\Gamma$ be a discrete group that acts continuously on $X$ from the right,
define $\widetilde{X} = \{(x,\gamma) \in X \times \Gamma : x\cdot\gamma= x\}$, and let $\Gamma$ act on $\widetilde{X}$
via the formula $(x,\gamma)\cdot\alpha = (x\cdot\alpha, \alpha^{-1}\gamma\alpha)$.  Results of P. Baum and A. Connes, 
along with facts about the Chern character, imply that 
$K^i_\Gamma(X) \otimes \mathbb{C} \cong K^i(\widetilde{X}\slash\Gamma) \otimes \mathbb{C}$ for $i = 0, -1$.  
In this note, we present an example where the groups $K^i_\Gamma(X)$ and $K^i(\widetilde{X}\slash\Gamma)$
are not isomorphic.
\end{abstract}
\maketitle

Let $\Gamma$ be a finite discrete group acting continuously on a compact Hausdorff space $X$ from the right. 
Define $\widetilde{X} = \{(x, \gamma) \in X \times \Gamma : x\cdot \gamma = x\}$ and endow $\widetilde{X}$ with the subspace topology
that it inherits from $X \times \Gamma$.  The group $\Gamma$ acts on $\widetilde{X}$ via the formula
$(x, \gamma)\cdot\alpha = (x\cdot\alpha, \alpha^{-1}\gamma\alpha)$, and we can consider the orbit space
$\widetilde{X}\slash\Gamma$.   Theorem 1.19 in \cite{BC} states that there exist isomorphisms
\begin{align*}
K^0_\Gamma(X) \otimes \mathbb{C} &\cong \sum_{j=0}^\infty H^{2j}(\widetilde{X}\slash\Gamma; \mathbb{C}) \\
K^{-1}_\Gamma(X) \otimes \mathbb{C} &\cong \sum_{j=0}^\infty H^{2j+1}(\widetilde{X}\slash\Gamma; \mathbb{C}),
\end{align*}
where $H^*(\widetilde{X}\slash\Gamma; \mathbb{C})$ denotes the {\v C}ech cohomology of $\widetilde{X}\slash\Gamma$
with complex coefficients.   We also have the Chern character isomorphisms 
\begin{align*}
K^0(\widetilde{X}\slash\Gamma) \otimes \mathbb{C} &\cong \sum_{j=0}^\infty H^{2j}(\widetilde{X}\slash\Gamma; \mathbb{C}) \\
K^{-1}(\widetilde{X}\slash\Gamma) \otimes \mathbb{C} &\cong \sum_{j=0}^\infty H^{2j+1}(\widetilde{X}\slash\Gamma; \mathbb{C}).
\end{align*}
Therefore $K^i_\Gamma(X) \otimes \mathbb{C} \cong K^i(\widetilde{X}\slash\Gamma) \otimes \mathbb{C}$ for $i = 0, -1$.
In other words, the groups $K^i_\Gamma(X)$ and $K^i(\widetilde{X}\slash\Gamma)$
are isomorphic up to torsion.  In this note, we present an example where 
these groups do not have isomorphic torsion subgroups.

Our example is Example B in \cite{Paschke}.  Consider the unit $3$-sphere $S^3$ in $\mathbb{R}^4$ and define
$\alpha: S^3 \longrightarrow S^3$ by $\alpha(x, y, z, t) = (-x, - y, -z, t)$.  The map $\alpha$ defines a $\mathbb{Z}_2$ action on
$S^3$.  From \cite{Green}, we know that the $\Gamma$-equivariant $K$-theory groups of a compact Hausdorff space $X$ are isomorphic 
to the operator algebra $K$-theory groups of $C(X) \rtimes \Gamma$.  Combining this fact with the computations in \cite{Paschke} we have 
$K^0_{\mathbb{Z}_2}(S^3) \cong \mathbb{Z}^3$ and $K^{-1}_{\mathbb{Z}_2}(S^3) \cong 0$.

Observe that $\widetilde{S^3}$ is homeomorphic to the disjoint union of $S^3$ and the fixed point set $F$ of the action of $\alpha$, and hence
$\widetilde{S^3}\slash\mathbb{Z}_2$ is homeomorphic
to the disjoint union of $S^3\slash\mathbb{Z}_2$ and $F$.   In our example, $F$ is the two-point set $\{(0, 0, 0, 1), (0, 0, 0, -1)\}$,
so $K^0(F) \cong \mathbb{Z}^2$ and $K^{-1}(F) \cong 0$.

To compute the $K$-theory of $S^3\slash\mathbb{Z}_2$, define closed sets 
\begin{align*}
A &= \{(x, y, z, t) \in S^3 : t \geq 0\} \\
B &= \{(x, y, z, t) \in S^3 : t \leq 0\}.
\end{align*}
Then 
$(A\slash\mathbb{Z}_2) \cup (B\slash\mathbb{Z}_2) = S^3\slash\mathbb{Z}_2$ and 
$(A\slash\mathbb{Z}_2) \cap (B\slash\mathbb{Z}_2) \cong S^2\slash\mathbb{Z}_2 \cong \mathbb{R}P^2$.  
Applying the Mayer-Vietoris sequence for reduced $K$-theory (\cite{Park}, Exercise 3.2), we have the six-term
exact sequence
\[
\xymatrix@R=30pt@C=30pt
{\widetilde{K}^0(S^3\slash\mathbb{Z}_2) \ar[r]^{} \ar[d]^{} 
&\widetilde{K}^0(A\slash\mathbb{Z}_2) \oplus \widetilde{K}^0(B\slash\mathbb{Z}_2) \ar[r]^{} 
&\widetilde{K}^0(\mathbb{R}P^2) \ar[d]^{} \\
\widetilde{K}^{-1}(\mathbb{R}P^2) 
&\widetilde{K}^{-1}(A\slash\mathbb{Z}_2) \oplus \widetilde{K}^{-1}(B\slash\mathbb{Z}_2) \ar[l]^{}
& \widetilde{K}^{-1}(S^3\slash\mathbb{Z}_2) \ar[l]^{}.}
\]
Both $A\slash\mathbb{Z}_2$ and $B\slash\mathbb{Z}_2$ are homeomorphic to the cone over $\mathbb{R}P^2$,
which has trivial reduced $K$-theory groups, and so the vertical maps in the six-term exact sequence are isomorphisms.  Therefore
\[
K^0(S^3\slash\mathbb{Z}_2) \cong \widetilde{K}^0(S^3\slash\mathbb{Z}_2) \oplus \mathbb{Z}
\cong \widetilde{K}^{-1}(\mathbb{R}P^2) \oplus \mathbb{Z} \cong  \mathbb{Z}
\]
\[
K^{-1}(S^3\slash\mathbb{Z}_2) \cong \widetilde{K}^{-1}(S^3\slash\mathbb{Z}_2) \cong 
\widetilde{K}^0(\mathbb{R}P^2) \cong \mathbb{Z}_2,
\]
which yield
\[
K^0(\widetilde{S^3}\slash\mathbb{Z}_2) \cong K^0(S^3\slash\mathbb{Z}_2) \oplus K^0(F) 
\cong \mathbb{Z} \oplus \mathbb{Z} \oplus \mathbb{Z} 
\]
\[
K^{-1}(\widetilde{S^3}\slash\mathbb{Z}_2) \cong K^{-1}(S^3\slash\mathbb{Z}_2) \oplus K^{-1}(F) \cong \mathbb{Z}_2.
\]
Thus $K^{-1}_{\mathbb{Z}_2}(S^3) \otimes \mathbb{C} \cong K^{-1}(\widetilde{S^3}\slash\mathbb{Z}_2) \otimes \mathbb{C}$,
but the groups $K^{-1}_{\mathbb{Z}_2}(S^3)$ and $K^{-1}(\widetilde{S^3}\slash\mathbb{Z}_2)$ are not isomorphic.

\end{document}